\documentclass[12pt,a4paper]{article}

\usepackage{amsmath,amsfonts}
\usepackage{mathrsfs} 
\usepackage{parskip}
\usepackage[T1]{fontenc}
\usepackage[utf8]{inputenc}
\usepackage{authblk}
\usepackage{enumerate}
\usepackage{graphicx}
\usepackage{xcolor}

\newtheorem{theorem}{Theorem}

\newtheorem{proposition}{Proposition}

\let\scr\mathscr

\def\Pb{\mathbf{P}}

\def\Ex{\mathbf{E}}
\def\Pb{\mathbf{P}}

\def\SS{\mathbb{S}}

\def\RR{\mathbb{R}}
\def\SS{\mathbb{S}}

\def\QQ{\mathbb{Q}}

\def\1{\mbox{1\hspace{-.25em}I}}
\begin{document}
\title{On Parameter Estimation of Hidden Ergodic Ornstein-Uhlenbeck Process}
 \author{ \textsc{Yury A. Kutoyants}\\ {\small Le Mans University, Le Mans, France }\\ 
 }

\date{}

\maketitle
\begin{abstract}
We consider the problem of parameter estimation for the partially observed 
linear stochastic differential equation. We assume that the unobserved
Ornstein-Uhlenbeck process depends on some unknown parameter and  
estimate the unobserved process and the unknown parameter simultaneously.  
We construct the two-step MLE-process for the estimator of the parameter and describe
its large sample asymptotic properties, including consistency and asymptotic normality. 
Using the Kalman-Bucy filtering equations we construct recurrent estimators of the state and the parameter. 
\end{abstract}


\section{Introduction}
We are given a partially observed linear system, defined by the equations
\begin{align}
\label{1}
{\rm d}X_t&=a Y_t\,{\rm d}t+
\sigma\, {\rm d}W_t,\qquad\qquad X_0=0, \\ {\rm
  d}Y_t&=-f Y_t\,{\rm d}t+b\,
           {\rm d}V_t,\qquad\qquad Y_0=\xi ,
\label{2}
\end{align}
 where $a\not=0, \sigma\not=0,b\not=0 $ and $f>0$ are   constants,
 $W^T=\left(W_t, 0\leq t\leq T\right)$ and $V^T=\left(V_t, 0\leq t\leq T\right)$
 are two independent Wiener 
 processes. The random variable  $\xi\sim{\cal N} \left(0, d^2\right)$ is
 independent of $W^T$ and $V^T$.

The system \eqref{1}-\eqref{2} is defined by the four parameters $a,f ,b,\sigma
^2$. Recall that the parameter $\sigma ^2$ can be estimated without error by
continuous time observations $X^T$ as follows. By the It\^o formula we can write
\begin{align*}
X_t^2&=2\int_{0}^{t}X_s\,{\rm d}X_s+\sigma ^2t
\end{align*}
Hence, for any $t\in (0,T]$, we have the estimator
\begin{align*}
\hat \sigma_t ^2=t^{-1}  X_t^2-2t^{-1}\int_{0}^{t}X_s\,{\rm d}X_s=\sigma ^2
\end{align*}
and this estimator equals the true value. Therefore we  consider only the
estimation of the three other parameters $f, b$ and $a$. Note that the
consistent estimation of the three-dimensional parameter $\vartheta
=\left(a,b,f\right)$ is impossible because the observed process can be written
as follows
\begin{align*}
X_t=ab\int_{0}^{t}\int_{0}^{s} e^{-f\left(s-r\right)}{\rm d}V_r\,{\rm
  d}s+\sigma W_t+o\left(1\right). 
\end{align*}
This means that the parameters $a$ and $b$ appear as product $ab$. We can have
consistent estimation of two-dimensional parameters $\left(f,a\right)$ and
$\left(f,b\right)$. This possibility we discuss in the last section.

The observations are $X^T=\left(X_t,0\leq t\leq T\right)$ and the
Ornstein-Uhlenbeck process $Y^T$ is unobservable (hidden), i.e., we have
partially observed linear model of observations.

 We consider estimation of the one-dimensional parameters $f$, $b$ and
 $a$ separately given the continuous time observations
 $X^T$.  The unknown parameter will be denoted by $\vartheta$ and we will assume that $\vartheta \in \Theta =\left(\alpha ,\beta
 \right)$ for some constants $\alpha<\beta$.  In all the cases the set $\Theta $ does not contain 0. Thus we are faced with 
 three different problems: $\vartheta =f$, $\vartheta =b$ and $\vartheta
 =a$. In each problem we propose a two-step construction of asymptotically
 efficient estimator-process of recurrent nature. First we propose a
 preliminary consistent estimator $\vartheta _{T^\delta }$ based on the
 observations $X^{T^\delta }=\left(X_t,0\leq t\leq T^\delta \right)$ with
 $\delta \in \left(1/2, 1\right)$. Then this estimator is used for
 construction of One-step MLE-process, which has recurrent structure.
In the last section we discuss the possibilities of the joint estimation of
two dimensional parameters $\vartheta =\left(f,b\right)$ and $\vartheta
=\left(f,a\right)$.

  Equations \eqref{1}-\eqref{2} is a prototypical model in the Kalman-Bucy
  filtering theory, which provides a closed form system of equations for the
  conditional expectation $m\left(t\right)=\Ex \left(Y_t |X_s,0\leq s\leq
  t\right)$ (\cite{A83}, \cite{K-B61},\cite{LS}). The statistical problems for
  discretely observed hidden Markov processes were studied by many authors
  (see \cite{BRR98}, \cite{CMT05}, \cite{EAM95}, \cite{EM02} and the
  references therein).  However, the literature on continuous time models is
  limited.  For the results in continuous time setup, we refer the interested
  reader to \cite{Kut94} (linear and non linear partially observed systems
  with small noise), \cite{EAM95} (continuous-time hidden Markov models
  estimation), \cite{PCh09} and \cite{KhK18} (hidden telegraph process
  observed in the white Gaussian noise).

In the present paper we are particularly interested in the asymptotic behavior
of the maximum likelihood estimator (MLE) $\hat\vartheta _T$ in the {\it large
  sample} asymptotic regime, i.e., when $T \rightarrow \infty $.  The
statistical problems for such observation models have been widely studied,
motivated by the importance of the Kalman-Bucy filtering in engineering
applications.

 Let us now recall the definitions of the MLE in the case $\vartheta =f$, when
 the other two parameters $a$ and $b$ are known. As the parameters of the
 model take finite values and $\sigma ^2>0$, the measures
 $\left\{\Pb_\vartheta^{\left(T\right)} ,\vartheta \in\Theta \right\}$ induced
 by the observations \eqref{1} on the space of continuous functions on
 $\left[0,T\right]$ are equivalent.  The likelihood ratio function (\cite{LS})
 is given by the expression
\begin{align}
\label{3}
L\left(\vartheta ,X^T\right)&=\exp\left\{\int_{0}^{T} \frac{a\,
  m\left(\vartheta,t \right)}{\sigma 
  ^2}\;{\rm d}X_t -\int_{0}^{T} \frac{a^2m\left(\vartheta,t
  \right)^2}{2 \sigma^2}\;{\rm d}t \right\}, 
\quad \vartheta \in \Theta .
\end{align}
Then the MLE $\hat\vartheta _T $ is defined by the equation
\begin{align}
\label{4}
L(\hat\vartheta _T ,X^T)=\sup_{\vartheta \in\Theta }L\left(\vartheta
,X^T\right) .
\end{align}
This means that to calculate $\hat\vartheta _T$ we need the values of the
family of stochastic processes $\left(m\left(\vartheta ,t\right), 0\leq t\leq
T\right),\vartheta \in\Theta $. The random process $m\left(\vartheta ,\cdot
\right)$ is solution of the Kalman-Bucy filtering equations (see \cite{A83},
\cite{K-B61}, \cite{LS})
\begin{align}
\label{5}
{\rm d}m\left(\vartheta ,t\right)&=-\vartheta \,m\left(\vartheta ,t\right){\rm
  d}t + \frac{\gamma \left(\vartheta ,t\right)a}{\sigma^2} \left[{\rm
    d}X_t-a\;m\left(\vartheta ,t\right){\rm d}t\right] ,\nonumber\\ &
=-\left[\vartheta +\frac{ \gamma \left(\vartheta ,t\right)a^2
  }{\sigma^2}\right] m\left(\vartheta ,t\right){\rm d}t + \frac{\gamma
  \left(\vartheta ,t\right)a}{\sigma^2} {\rm d}X_t,
\end{align}
where $m\left(\vartheta ,0\right)= \Ex_\vartheta \left(\xi |X_0\right)=0$. The
function $\gamma \left(\vartheta ,t\right)=\Ex_\vartheta
\left(m\left(\vartheta ,t\right)-Y_t \right)^2 $ is the solution of the
Ricatti equation
\begin{align}
\label{6}
\frac{\partial \gamma \left(\vartheta ,t\right)}{\partial t}=
-2\vartheta\,\gamma \left(\vartheta ,t\right) -\frac{\gamma \left(\vartheta
  ,t\right)^2a^2}{\sigma^2} +b^2,\qquad \gamma \left(\vartheta ,0\right)=d^2 .
\end{align}
Due to importance of this model in many applied problems, much
engineering literature is concerned with identification of this model.

The  behavior of the MLE was studied at least in three asymptotics:
\begin{itemize}
\item  {\it Small noise in
both equations } $\sigma =b=\varepsilon\rightarrow 0 $ ($T$ is fixed)
\cite{Kut84}, \cite{Kut94} 
\begin{align*}
\frac {\hat\vartheta _\varepsilon -\vartheta }{\varepsilon }\Longrightarrow
      {\cal N}\left(0, {\rm I}\left(\vartheta \right)^{-1}\right).
 \end{align*}

\item {\it Large sample} $T\rightarrow \infty $ ($\sigma $ and $b$ are fixed)
  \cite{Kut04}
\begin{align*}
\sqrt{T} \left({\hat\vartheta _T -\vartheta }\right)\Longrightarrow
      {\cal N}\left(0, {\rm I}\left(\vartheta \right)^{-1}\right).
\end{align*}

\item {\it Small noise in observation} only, $\sigma \rightarrow 0$, ($T$ and
  $b$ are fixed) \cite{K19}
\begin{align*}
\frac {\hat\vartheta _\varepsilon -\vartheta }{\sqrt{\varepsilon}
}\Longrightarrow {\cal N}\left(0, {\rm I}\left(\vartheta \right)^{-1}\right).
 \end{align*} 
\end{itemize}
In all three cases ${\rm I}\left(\vartheta \right) $ the Fisher information is
different. It was also shown that the polynomial moments of the scaled
estimation error converge and the MLE is asymptotically efficient.

It is evident that the numerical calculation of the MLE $\hat\vartheta _T$
according to \eqref{3}-\eqref{6} is quite a difficult problem. The goal of
this work is to suggest the new estimator, called One-step MLE-process
$\vartheta _t^\star , \tau \leq t\leq T$, which has two advantages. First, its
numerical calculation is much more simple than that of the MLE and, second,
this estimator has a recurrent structure and can be used for the joint
estimation of the hidden process $Y_t$ and the parameter $\vartheta $. Similar
One-step MLE's and Multi-step MLE-processes, introduced in \cite{Kut17}, have
been applied in the problem of parameter estimation of the hidden telegraph
process \cite{KhK18}, parameter estimation in diffusion processes by the
discrete time observations \cite{KU15}, in the problem of frequency estimation
\cite{G84}, intensity parameter estimation for inhomogeneous Poisson processes
\cite{DGK17}, parameter estimation for the Markov sequences \cite{KM16}.
 
\section{Preliminary estimator.}

Following \cite{KhK18} One-step MLE process will be constructed in two
steps. First we introduce a 
consistent and asymptotically normal preliminary  estimator and then this
estimator is used to define One-step MLE-process. Preliminary estimator is
constructed using an asymptotically negligible amount of the observations 
$X^{K } =\left(X_t,0\leq t\leq K \right)$, where $K =T^\delta ,
\delta \in \left(1/2,1\right)$.

Suppose that $\vartheta =f$ and introduce the statistic $\SS_K$ and the function
$\Phi \left(\vartheta \right),\vartheta \in\Theta $:
\begin{align*}
\SS_K&=\frac{1}{K}\sum_{k=1}^{K}\left[X_k-X_{k-1}\right]^2,\qquad \quad \Phi
\left(\vartheta \right)=\frac{a^2b^2}{\vartheta ^3}\left[e^{-\vartheta
  }-1+\vartheta \right]+\sigma ^2.
\end{align*}
In the cases $\vartheta =b$ and $\vartheta =a$ the counterparts of the latter function are
\begin{align*}
\Phi_*
\left(\vartheta \right)=\frac{a^2\vartheta ^2}{f ^3}\left[e^{-f
  }-1+f \right]+\sigma ^2,\qquad \hat\Phi
\left(\vartheta \right)=\frac{\vartheta ^2b^2}{f^3}\left[e^{-f
  }-1+f \right]+\sigma ^2
\end{align*}
respectively. 

In this section we consider the case $\vartheta =f$ only. Therefore 
\begin{align}
\label{7}
{\rm d}X_t&=a Y_t\,{\rm d}t+ \sigma\, {\rm d}W_t,\qquad\qquad X_0=0, \\
 {\rm   d}Y_t&=-\vartheta  Y_t\,{\rm d}t+b\, {\rm d}V_t,\qquad\qquad Y_0=\xi ,
\label{8}
\end{align}

Note that the function $\Phi \left(\vartheta \right), \alpha <\vartheta <\beta
$ is strictly decreasing. Define the preliminary estimator $\bar\vartheta _K$, 
base the observations $X^K $:
\begin{align*}
\bar\vartheta _K=\vartheta _K^* \1_{\left\{{\cal
    A}_K\right\}}+\alpha\1_{\left\{{\cal A}_K^-\right\}} +\beta
\1_{\left\{{\cal A}_K^+\right\}}.
\end{align*}
Here $\vartheta _K^*$ is the root of  equation   $\Phi \left(\vartheta
_K^*\right)=\SS_K$ and  ${\cal A}_K,{\cal A}_K^-,{\cal A}_K^+ $ are   the sets
\begin{align*}
{\cal A}_K&=\left\{\omega :\Phi \left(\beta\right) <\SS_K<\Phi \left(\alpha
\right) \right\},\quad {\cal A}_K^-=\left\{\omega : \SS_K\geq \Phi
\left(\alpha\right) \right\},\\ {\cal A}_K^+&=\left\{\omega :\SS_K\leq \Phi
\left(\beta \right) \right\} .
\end{align*}
The asymptotic behavior of $\bar\vartheta _K$ as $K\rightarrow \infty $ is
described in the following proposition. 
\begin{proposition}
\label{P1}
The estimator $\bar\vartheta _K$ is consistent, uniformly on compacts\\
$\left[\bar \alpha,\bar\beta \right]$ $\subset \Theta$, and
\begin{align}
\label{9}
\sup_{\vartheta _0\in\Theta } \Ex_{\vartheta _0}\left| \bar\vartheta
  _K-\vartheta _0\right|^2\leq \frac{C}{K}
\end{align}
with some constant $C>0$.
\end{proposition}
{\bf Proof}. We have
\begin{align*}
\Ex_{\vartheta _0}\left[ \bar\vartheta
  _K-\vartheta _0\right]^2&=\Ex_{\vartheta _0}\left[ \vartheta
  _K^*-\vartheta _0\right]^2\1_{\left\{{\cal
    A}_K\right\}}+\left(\vartheta _0-\alpha  \right)^2\Pb_{\vartheta
  _0}\left({\cal A}_K^- \right)\\
&\qquad +\left( \beta -\vartheta _0 \right)^2\Pb_{\vartheta _0}\left({\cal
  A}_K^+ \right) .
\end{align*}
For the probabilities we have the estimates
\begin{align*}
\Pb_{\vartheta _0}\left({\cal A}_K^- \right)&=\Pb_{\vartheta _0}\left(
\SS_K-\Phi \left(\vartheta _0\right) \geq \Phi \left(\alpha\right)-\Phi
\left(\vartheta _0\right) \right) \leq \frac{\Ex_{\vartheta
    _0}\left[\SS_K-\Phi \left(\vartheta _0\right) \right]^2 }{\left|\Phi
  \left(\alpha\right)-\Phi \left(\vartheta _0\right)\right|^2}
,\\ \Pb_{\vartheta _0}\left({\cal A}_K^+ \right)&\leq \frac{\Ex_{\vartheta
    _0}\left|\SS_K-\Phi \left(\vartheta _0\right) \right|^2 }{\left[\Phi
    \left(\beta \right)-\Phi \left(\vartheta _0\right)\right]^2}.
\end{align*}

Therefore we have to  study the asymptotics  
of the statistic $\SS_K$ as $K\rightarrow \infty $:
\begin{align*}
\SS_K&=\frac{1}{K}\sum_{k=1}^{K}\left[X_k-X_{k-1}\right]^2
=\frac{1}{K}\sum_{k=1}^{K}\left[\int_{k-1}^{k}{\rm
    d}X_s\right]^2\\
&=\frac{a^2}{K}\sum_{k=1}^{K}\eta _k^2 +\frac{2a\sigma
}{K}\sum_{k=1}^{K}\eta _k \left[W_k-W_{k-1}\right]+\frac{\sigma^2
}{K}\sum_{k=1}^{K} \left[W_k-W_{k-1}\right]^2,
\end{align*}
where
\begin{align*}
\eta _k =\int_{k-1}^{k}Y_t\;{\rm d}t.
\end{align*}
We have
\begin{align*}
\Ex_{\vartheta _0}\SS_K&=\frac{a^2}{K}\sum_{k=1}^{K}\Ex_{\vartheta _0}\eta _k^2+\sigma ^2
\end{align*}
because $Y^T=\left(Y_t,0\leq t\leq T\right)$ and $W^T=\left(W_t,0\leq t\leq
T\right)$ are independent. 

The process $Y^T$ can be written as 
\begin{align*}
Y_t=\xi e^{-{\vartheta_0} t}+b\int_{0}^{t}e^{-{\vartheta_0} \left(t-r\right)}{\rm d}V_r.
\end{align*}
Hence
\begin{align*}
\Ex_{\vartheta _0}Y_tY_s&=\Ex_{\vartheta _0}\xi ^2e^{-{\vartheta_0}
  \left(t+s\right)}+b^2e^{-{\vartheta_0} \left(t+s\right)}
\int_{0}^{t\wedge s}e^{2{\vartheta_0} r}{\rm d}r\\
&=\left[d^2-\frac{b^2}{2\vartheta _0} \right] e^{-{\vartheta_0}
  \left(t+s\right)}+\frac{b^2}{2\vartheta _0}e^{-{\vartheta_0}
  \left|t-s\right|}
\end{align*}
and
\begin{align*}
\Ex_{\vartheta _0}\eta _k^2& =\int_{k-1}^{k}\int_{k-1}^{k}\Ex_{\vartheta
  _0}Y_tY_s\,{\rm d}s{\rm d}t \\ &=\left[d^2-\frac{b^2}{2\vartheta _0}
  \right]\left(\int_{k-1}^{k}e^{-{\vartheta_0} t}{\rm d}t\right)^2
+\frac{b^2}{2\vartheta _0}\int_{k-1}^{k}\int_{k-1}^{k}e^{-{\vartheta_0}
  \left|t-s\right|}{\rm d}s{\rm d}t\\ &=\left[\frac{d^2}{\vartheta
    _0^2}-\frac{b^2}{2\vartheta _0^3} \right] \left[e^{\vartheta_0}
  -1\right]^2e^{-2{\vartheta_0} k}+\frac{b^2}{{\vartheta_0}
  ^3}\left[e^{-{\vartheta_0} }-1+{\vartheta_0} \right].
\end{align*}
Therefore
\begin{align*}
\Ex_{\vartheta _0}\SS_K&=\left[e^{\vartheta_0} -1\right]^2\frac{a^2d^2}{{\vartheta_0}
  ^2K}\sum_{k=1}^{K}e^{-2{\vartheta_0} k}+\frac{a^2b^2}{{\vartheta_0}
  ^3}\left[e^{-{\vartheta_0} }-1+{\vartheta_0} \right]+\sigma ^2 \\
&=\frac{a^2b^2}{{\vartheta_0}
  ^3}\left[e^{-{\vartheta_0} }-1+{\vartheta_0} \right]+\sigma ^2+r_K ,\qquad
\left|r_K\right|\leq \frac{C}{K}. 
\end{align*}
Using similar calculations we obtain the estimate
\begin{align*}
\Ex_{\vartheta _0}\eta _k\eta _m\leq C\,e^{-{\vartheta_0} \left|k-m\right|}
\end{align*}
which allows us to prove the law of large numbers: for $K\rightarrow
\infty $ we have convergence in mean square 
\begin{align*}
\Ex_{\vartheta _0}\left(\SS_K-\Ex_{\vartheta _0}\SS_K \right)^2\leq
\frac{C}{K},\qquad \SS_K\longrightarrow \Phi \left({\vartheta_0}
\right)=\frac{a^2b^2}{{\vartheta_0} ^3}\left[e^{-{\vartheta_0}
  }-1+{\vartheta_0} \right]+\sigma ^2
\end{align*}
and
\begin{align*}
\Ex_{\vartheta _0}\left(\SS_K-\Phi \left({\vartheta_0} \right)\right)^2\leq
2\Ex_{\vartheta _0}\left(\SS_K-\Ex_{\vartheta _0}\SS_K \right)^2 +2
\left(\Ex_{\vartheta _0}\SS_K-\Phi \left({\vartheta_0} \right)\right)^2\leq
\frac{C}{K}.
\end{align*}
Hence
\begin{align*}
\sup_{\bar\alpha <\vartheta _0\leq \bar\beta }\left[\Pb_{\vartheta
    _0}\left({\cal A}_K^- \right)+\Pb_{\vartheta _0}\left({\cal A}_K^+
  \right)\right]\leq \frac{C}{K}.
\end{align*}
The function $\Phi \left(\vartheta \right),\alpha <\vartheta <\beta $ is
strictly decreasing. If we denote its inverse function as $\Psi
\left(\phi\right)=\Phi ^{-1}\left(\phi \right),\Phi \left(\beta
\right)<\phi<\Phi\left(\alpha \right) $, then we have
\begin{align*}
\Psi'\left(\Phi\right) =\frac{1}{\Phi '\left(\vartheta \right)}, \qquad {\rm
  for} \qquad \Phi =\Phi \left(\vartheta \right)
\end{align*}
and 
\begin{align*}
\sup_{\vartheta \in\Theta } \left|\Psi'\left(\Phi\left(\vartheta
\right)\right)\right|=\left(\inf_{\vartheta \in\Theta }\left|\Phi
'\left(\vartheta \right) \right|\right)^{-1}= \left|\Phi
'\left(\beta  \right)\right|^{-1}\equiv c_*>0.
\end{align*}

We can write
\begin{align*}
\Ex_{\vartheta _0}\left[ \vartheta _K^*-\vartheta _0\right]^2\1_{\left\{{\cal
    A}_K\right\}}&=\Ex_{\vartheta _0}\left[ \Psi \left( \SS_K \right)-\Psi
  \left(\Phi \left(\vartheta _0\right)\right)\right]^2\1_{\left\{{\cal
    A}_K\right\}}\\
&\leq c_*^{-2} \Ex_{\vartheta _0}\left[\SS_K -\Phi \left(\vartheta _0\right)
  \right]^2\leq \frac{C}{c_*^2K}\longrightarrow 0 
\end{align*}
as $K\rightarrow \infty $.

If we put $K=T^\delta $, then
\begin{align}
\label{10}
\sup_{\bar\alpha \leq \vartheta _0<\bar\beta }\Ex_{\vartheta _0}\left[
  \vartheta _{T^\delta }^*-\vartheta _0\right]^2\leq C\,T^{-\delta }.
\end{align}

\section{One-Step MLE-process. Case $\vartheta =f$.}

Suppose that the unknown parameter is $\vartheta =f$ and we have the model
\eqref{7}-\eqref{8}, where the process $X^T$ is observable and the
Ornstein-Uhlenbeck process $Y^T$ is ``hidden''.  We realize the asymptotically
efficient estimation of the parameter $\vartheta \in\Theta $ in two
steps. First we calculate the preliminary estimator $\bar\vartheta _{T^\delta
} $ and then using this estimator we construct the One-step MLE-process.

Recall that the equation \eqref{6} has explicit solution
\begin{align*}
\gamma \left(\vartheta ,t\right)=e^{-2r\left(\vartheta
  \right)t}\left[\frac{1}{\gamma _0-\gamma \left(\vartheta \right)} +
  \frac{a^2}{2r\left(\vartheta \right)\sigma ^2}\left(1-e^{-2r\left(\vartheta
    \right)t}\right)\right]^{-1}+ \gamma \left(\vartheta \right).
\end{align*}
Here $\gamma _0=d^2$,
\begin{align*}
r\left(\vartheta \right)=\left( \vartheta ^2+\frac{b^2a^2}{\sigma
  ^2}\right)^{1/2},\qquad  \gamma \left(\vartheta \right)=\frac{\vartheta
  \sigma ^2}{a^2}\left(\sqrt{1+\frac{b^2a^2}{\vartheta ^2\sigma ^2}}-1\right). 
\end{align*}
Therefore we have exponential convergence of $\gamma \left(\vartheta
,t\right)$ to the stationary solution $\gamma\left(\vartheta\right)$
\begin{align*}
\left|\gamma \left(\vartheta ,t\right)-\gamma \left(\vartheta
\right)\right|\leq C\,e^{-2r\left(\vartheta \right)t}. 
\end{align*}
To simplify the exposition we suppose that $d^2=\gamma \left(\vartheta
\right)$; then we have $\gamma \left(\vartheta ,t\right)=\gamma \left(\vartheta
\right) $. The case with an arbitrary $d^2$ requires cumbersome calculations,
but  the main results remain intact.  

The equation for $m\left(\vartheta ,t\right)$ in this case is 
\begin{align*}
{\rm d}m\left(\vartheta ,t\right)&=-\left[\vartheta +\frac{ \gamma
    \left(\vartheta \right)a^2 }{\sigma^2}\right] m\left(\vartheta
,t\right){\rm d}t + \frac{\gamma \left(\vartheta\right)a}{\sigma^2} {\rm
  d}X_t,
\end{align*}

Denote $m_t=m\left(\vartheta _0,t\right)$ 
and $\gamma _*\left(\vartheta_0\right)=\gamma _*$, 
where $\vartheta _0$ is the true value. Then for the
process $m_t,0\leq t\leq T$ we obtain the equation
\begin{align}
\label{11}
{\rm d}m_t=-\vartheta_0 m_t{\rm d}t+\frac{\gamma_*a}{\sigma }{\rm d}\bar
W_t,\quad m_0\sim{\cal N}\left(0,\gamma_* \right),\quad 0\leq t\leq T,
\end{align}
Here 
we used the {\it innovation theorem} (see \cite{LS}, ??)
\begin{align*}
{\rm d}X_t=am_t\,{\rm d}t+\sigma\, {\rm d}\bar W_t,\qquad X_0=0, \qquad 0\leq
t\leq T.
\end{align*}
The {\it innovation} Wiener process $\bar W_t$ is defined by this equation and
$m_0$ is independent on $\bar W_t,0\leq t\leq T$.  With probability 1, the
random process $m\left(\vartheta ,t\right)$ has continuous derivatives
w.r.t. $\vartheta $ and derivative processes $\dot m\left(\vartheta
,t\right),\ddot m\left(\vartheta ,t\right) $ satisfy the equations
\begin{align}
\label{12}
{\rm d}\dot m\left(\vartheta ,t\right)&=-\left[\vartheta+\frac{ \gamma
    \left(\vartheta \right)a^2}{\sigma^2} \right]\,\dot m\left(\vartheta
,t\right){\rm d}t + \frac{\dot \gamma \left(\vartheta \right)a}{\sigma^2} {\rm
  d}X_t \nonumber\\ &\qquad  -\left[ 1+ \frac{\dot \gamma \left(\vartheta
    \right)a^2}{\sigma^2} \right] m\left(\vartheta ,t\right){\rm d}t,\\
\label{13}
{\rm d}\ddot m\left(\vartheta ,t\right)&=-\left[\vartheta+\frac{ \gamma
    \left(\vartheta \right)a^2}{\sigma^2} \right]\,\ddot m\left(\vartheta
,t\right){\rm d}t + \frac{\ddot \gamma \left(\vartheta \right)a}{\sigma^2}
{\rm d}X_t \nonumber\\ &\qquad -2\left[ 1+ \frac{\dot \gamma \left(\vartheta
    \right)a^2}{\sigma^2} \right] \dot m\left(\vartheta ,t\right){\rm
  d}t-\frac{\ddot \gamma \left(\vartheta \right)a^2}{\sigma^2}
m\left(\vartheta ,t\right){\rm d}t,
\end{align}
The Fisher information for this model of observations is
\begin{align*}
{\rm I}\left(\vartheta \right)=\frac{1}{2\vartheta }-\frac{2\dot
  r\left(\vartheta \right)}{r\left(\vartheta \right)+\vartheta }+ \frac{\dot
  r\left(\vartheta \right)^2}{2r\left(\vartheta \right)}.
\end{align*}
Note that ${\rm I}\left(\vartheta \right)$ has continuous bounded derivatives
and is uniformly in $\vartheta \in \Theta $ separated from zero.

According to \cite{Kut17} the One-step MLE-process $\vartheta _t^\star,
T^\delta < t\leq T$ is introduced as follows
\begin{align}
\label{14}
\vartheta _t^\star=\bar\vartheta _{T^\delta }+\frac{a}{ \sigma ^2
  t {\rm I}\left(\bar\vartheta _{T^\delta } \right)}
\int_{T^\delta }^{t}\dot m(\bar\vartheta _{T^\delta },s) \left[{\rm
    d}X_s-am(\bar\vartheta _{T^\delta },s) {\rm d}s\right].
\end{align}

Let us  change the variables $t=\tau T$ and denote $\vartheta _{\tau
  T}^\star=\vartheta_T ^\star\left(\tau \right), T^{\delta -1}<\tau \leq 1 $. 
\begin{theorem}
\label{T1}  One-step MLE-process  $\vartheta_T ^\star\left(\tau \right),
T^{\delta -1}<\tau \leq 1  $ with $\delta \in \left(1/2, 1\right)$ is
consistent: for any $\nu >0$ and any $\tau \in (0,1]$
\begin{align*}
\lim_{T\rightarrow \infty }-\Pb_{\vartheta _0}\left(\left|\vartheta_T
^\star\left(\tau \right)-\vartheta _0 
\right|>\nu \right)=  0, 
\end{align*}
and asymptotically normal
\begin{align*}
\sqrt{\tau T}\left(\vartheta_T ^\star\left(\tau \right)-\vartheta _0
\right)\Longrightarrow {\cal N}\left(0, {\rm I}\left(\vartheta_0 \right)^{-
  1}\right).
\end{align*}
\end{theorem}
{\bf Proof.} Consider the difference
\begin{align}
\label{15}
&\sqrt{\tau T}\left(\vartheta_T ^\star\left(\tau \right)-\vartheta _0\right)=
\sqrt{\tau T}\left(\bar\vartheta _{T^\delta }-\vartheta _0\right)\nonumber\\ 
&\qquad \quad+\frac{a}{ \sigma \sqrt{\tau T} {\rm
    I}\left(\bar\vartheta _{T^\delta } \right)} \int_{T^\delta }^{\tau T}\dot
m(\bar\vartheta _{T^\delta },s) {\rm d}\bar
W_s\nonumber\\
 &\qquad \quad+\frac{a^2}{ \sigma ^2 \sqrt{\tau T} {\rm
    I}\left(\bar\vartheta _{T^\delta } \right)} \int_{T^\delta }^{\tau T}\dot
m\left(\bar\vartheta _{T^\delta },s\right) \left[m_s-m(\bar\vartheta
  _{T^\delta },s)\right] {\rm d}s.
\end{align}
Note that as it follows from the equations \eqref{12}-\eqref{13}, the Gaussian
processes $ \dot m\left(\vartheta ,t\right) $ and $ \ddot m\left(\vartheta
_0,t\right) $  have bounded variances and therefore for any   $p>1 $ we have 
\begin{align*}
\sup_{\vartheta \in\Theta }\Ex_{\vartheta _0} \left| \dot m\left(\vartheta
,t\right)\right|^p\leq C,\qquad \sup_{\vartheta \in\Theta }\Ex_{\vartheta _0}
\left| \ddot m\left(\vartheta ,t\right)\right|^p\leq C,
\end{align*}
where the constants do not depend on $t$. 
We can write
\begin{align*}
\dot m(\bar\vartheta _{T^\delta },s)&=\dot m(\vartheta _0,s)+\dot
m(\bar\vartheta _{T^\delta },s)-\dot m\left(\vartheta _0 ,s\right)\\ &=\dot
m\left(\vartheta _0,s\right)+\left(\vartheta _0- \bar\vartheta _{T^\delta }
\right)\ddot m(\tilde \vartheta ,s)=\dot
m\left(\vartheta _0,s\right)+O\left(T^{-\delta /2}\right)
\end{align*}
  because
\begin{align*}
 \left( \Ex_{\vartheta _0} \left|\left(\vartheta _0- \bar\vartheta _{T^\delta }
\right)\ddot m(\tilde \vartheta ,s)  \right|\right)^2\leq \Ex_{\vartheta _0}
\left(\vartheta _0- \bar\vartheta _{T^\delta } 
\right) ^2\Ex_{\vartheta _0} \ddot m(\tilde \vartheta ,s)  ^2\leq \frac{C}{T^\delta }.
\end{align*}
Here   $ \left|\tilde \vartheta -\bar\vartheta _{T^\delta }\right|\leq \left|
\vartheta_0 -\bar\vartheta _{T^\delta }\right|$. 

 Further, for the Fisher information we have 
\begin{align*}
\left|\frac{1 }{{\rm I}\left(\bar\vartheta _{T^\delta } \right)}-\frac{1}{{\rm
    I}\left(\vartheta _0 \right)}\right|=\frac{\left|{\rm
    I}\left(\bar\vartheta _{T^\delta } \right)-{\rm I}\left(\vartheta _0
  \right)\right| }{{\rm I}\left(\bar\vartheta _{T^\delta } \right){\rm
    I}\left(\vartheta _0 \right) }\leq C\left|\bar\vartheta _{T^\delta }
-\vartheta _0 \right|=O\left(T^{-\delta /2}\right)
\end{align*}
This allows us to write 
\begin{align*}
\Delta _T&=\frac{a}{ \sigma \sqrt{\tau T} {\rm
    I}\left(\bar\vartheta _{T^\delta } \right)} \int_{T^\delta }^{\tau T}\dot
m(\bar\vartheta _{T^\delta },s)\; {\rm d}\bar W_s\\
&\quad \qquad =\frac{a}{\sigma  {\rm I}\left(\vartheta _0\right)\sqrt{\tau
    T-T^\delta} } \int_{T^\delta }^{\tau
  T}\dot m(\vartheta _0,s)\; {\rm d}\bar W_s \left(1+o\left(1\right)\right).
\end{align*}
By the law of large numbers 
\begin{align*}
\frac{a}{\sigma \tau T}\int_{T^\delta }^{\tau
  T}\dot m(\vartheta _0,s)^2\; {\rm d}s\longrightarrow  {\rm I}\left(\vartheta _0\right)
\end{align*}
and therefore by the central limit theorem
\begin{align*}
\frac{a}{\sigma  {\rm I}\left(\vartheta _0\right)\sqrt{\tau
    T} } \int_{T^\delta }^{\tau
  T}\dot m(\vartheta _0,s)\; {\rm d}\bar W_s \Longrightarrow  {\cal
  N}\left(0,{\rm I}\left(\vartheta _0\right)^{-1} \right) .
\end{align*}
The similar arguments allow us to write
\begin{align*}
&\int_{T^\delta }^{\tau T}\dot m\left(\bar\vartheta _{T^\delta },s\right)
\left[m\left(\vartheta _0,s\right)-m(\bar\vartheta _{T^\delta },s)\right] {\rm
  d}s\\
&\qquad \qquad =-(\bar\vartheta _{T^\delta }-\vartheta _0)\int_{T^\delta }^{\tau T}\dot
m\left(\bar\vartheta _{T^\delta },s\right) \dot m(\tilde\vartheta,s)
  {\rm d}s\\
&\qquad \qquad =-(\bar\vartheta _{T^\delta }-\vartheta _0)\int_{T^\delta }^{\tau T}\dot
m\left(\vartheta _0,s\right)^2 {\rm d}s \left(1+O\left(T^{-\delta /2}\right)\right) .
\end{align*}
Recall that as we have stationary regime $\Ex_{\vartheta _0} \dot
m(\vartheta_0,s)^2= \sigma ^2a^{-2}{\rm I}\left(\vartheta _0\right)
$. Therefore
\begin{align*}
&\frac{1}{{\tau T}}\int_{T^\delta }^{\tau T}\dot m\left(\vartheta
_0,s\right)^2 {\rm d}s-\sigma ^2a^{-2}{\rm I}\left(\vartheta
_0\right)=\frac{1}{\sqrt{\tau T }}A\left(\tau T\right)
\end{align*}
where the integral (see, e.g., Proposition 1.23 in \cite{Kut04})
\begin{align*}
A\left(\tau T\right)=\frac{1}{\sqrt{\tau T}}\int_{T^\delta }^{\tau
  T}\left[\dot m\left(\vartheta _0,s\right)^2-\Ex_{\vartheta _0}\dot
  m\left(\vartheta _0,s\right)^2\right] {\rm d}s\Longrightarrow {\cal
  N}\left(0,D\left(\vartheta _0\right)\right).
\end{align*}
Hence we obtained the representation
\begin{align*}
&\frac{a^2}{ \sigma ^2 \sqrt{\tau T} {\rm
    I}\left(\bar\vartheta _{T^\delta } \right)} \int_{T^\delta }^{\tau T}\dot
m\left(\bar\vartheta _{T^\delta },s\right) \left[m_s-m(\bar\vartheta
  _{T^\delta },s)\right] {\rm d}s\\
&\qquad = -\sqrt{\tau T}\left(\bar\vartheta _{T^\delta }-\vartheta
_0\right)\left(1+O\left(T^{-\delta /2}\right)\right). 
\end{align*}
Substitution of this relation into the initial representation \eqref{15}
yields the final expression
\begin{align*}
\sqrt{\tau T}\left(\vartheta ^\star\left(\tau \right)-\vartheta _0\right)&=
\Delta _T+\sqrt{\tau T}\left(\bar\vartheta _{T^\delta }-\vartheta
_0\right)O\left(T^{-\delta /2}\right)= \Delta _T+O\left(T^{\frac{1}{2}-\delta
}\right)\\
&\Longrightarrow  {\cal
  N}\left(0,{\rm I}\left(\vartheta _0\right)^{-1} \right),
\end{align*}
since $\delta \in (1/2,1)$.

Note that the process $\vartheta _t^\star, T^\delta <t\leq T$ can be
written in recurrent form
\begin{align}
\label{16}
{\rm d}\vartheta _t^\star=-\frac{\vartheta _t^\star}{t-T^\delta
}{\rm d}t +\frac{a\dot m(\bar\vartheta
  _{T^\delta },t)}{\sigma ^2t{\rm I}(\bar\vartheta
  _{T^\delta }) } \left[{\rm d}X_t-a m(\bar\vartheta
  _{T^\delta },t){\rm d}t\right] 
\end{align}
and we can introduce the adaptive filtering equations as follows 
\begin{align}
\label{17}
{\rm d}m_t&=-\left[\vartheta _t^\star +\frac{ \gamma \left(\vartheta _t^\star
    \right)a^2 }{\sigma^2}\right] m_t{\rm d}t +
\frac{\gamma \left(\vartheta _t^\star \right)a}{\sigma^2} {\rm d}X_t,\quad
T^\delta <t\leq T, \\
\gamma  \left(\vartheta _t^\star \right)&=\frac{\vartheta _t^\star\sigma
  ^2}{a^2}\left( \sqrt{1+\frac{b^2a^2}{\left(\vartheta _t^\star\right)^2\sigma
    ^2 }}-1\right) 
\label{18}
\end{align}
with the initial value $m_{T^\delta }=m\left(\bar\vartheta
  _{T^\delta },T^\delta \right) $. Here 
$$
m_t=\left.\Ex_\vartheta \left(Y_t|X_s,0\leq
s\leq t\leq t\right)\right|_{\vartheta =\bar\vartheta
  _{T^\delta }}.
$$  
It will be interesting to see the behavior of the system \eqref{16}-\eqref{18}
using   numerical simulations. 

Recall that if we put $\tau =1$, then $\vartheta _T^\star$  is One-step MLE
with 
\begin{align*}
\sqrt{T}\left(\vartheta _T^\star-\vartheta _0 \right)\Longrightarrow {\cal
  N}\left(0, {\rm I}\left(\vartheta _0\right)^{-1}\right)
\end{align*}
studied for ergodic diffusion processes in the Section 2.5 \cite{Kut04}.
Therefore the estimator $\vartheta _T^\star $ is asymptotically equivalent to
the asymptotically efficient MLE $\hat\vartheta _T$ defined by the equation
\eqref{4}. There is essential  computational difference between these two
estimators. The calculation of $\hat\vartheta _T$ using
\eqref{3}-\eqref{6} requires solving the differential equations
\eqref{5}-\eqref{6} for {\em numerous} values of $\vartheta \in\Theta $, which is computationally 
inefficient. 
 To construct One-step MLE-process $\vartheta _T^\star $ we have to calculate a simple preliminary
estimator $\bar\vartheta _{T^\delta }$ and then to solve the system
\eqref{5}-\eqref{6}  for just one value $\vartheta =\bar\vartheta _{T^\delta
}$. The difference between these two approaches becomes even more significant in
the case of multidimensional $\vartheta $.

\section{One-Step MLE-process. Case $\vartheta =b$.}

Suppose that the volatility $b=\vartheta $ is the unknown parameter and we have the equations
\begin{align}
\label{19}
{\rm d}X_t&=a Y_t\,{\rm d}t+
\sigma\, {\rm d}W_t,\qquad\qquad X_0=0, \\ {\rm
  d}Y_t&=-f Y_t\,{\rm d}t+\vartheta \,
           {\rm d}V_t,\qquad\qquad Y_0=\xi.
\label{20}
\end{align}
As before all parameters $a,\sigma,\vartheta $ do not vanish 
and $f>0$.  The volatility $\vartheta \in
\left(\alpha ,\beta \right)$ with $\alpha >0$ and the function 
\begin{align*}
\Phi_*\left(\vartheta \right)=\frac{a^2\vartheta
  ^2}{f^3}\left[e^{-f}-1+f\right] +\sigma ^2,\qquad \alpha <\vartheta <\beta 
\end{align*}
is strictly increasing.

The statistic $\SS_K$, with the new notations, converges to this function
\begin{align*}
\SS_K\longrightarrow \Phi_*\left(\vartheta_0 \right)\qquad {\rm as} \qquad
K\rightarrow \infty .
\end{align*}
Therefore we have the explicit expression for the preliminary estimator
\begin{align*}
\bar\vartheta_{K } =\vartheta _K^* \1_{\left\{{\cal B}_K\right\}
} +\alpha  \1_{\left\{{\cal B}_K^-\right\}
}+\beta \1_{\left\{{\cal B}_K^+\right\}
}      ,
\end{align*}
where
\begin{align*}
\vartheta _K^*=\left(\frac{f^3\left(\SS_K-\sigma
  ^2\right)}{a^2\left[e^{-f}-1+f\right] }\right)^{1/2}.
\end{align*}
Here the sets ${\cal B}^\pm$ are defined by the similar relations
\begin{align*}
 {\cal B}_K^-&=\left\{\omega :
\SS_K\leq \Phi _*\left(\alpha \right)\right\},\qquad {\cal B}_K^+=\left\{\omega :
\SS_K\geq \Phi _*\left(\beta  \right)\right\},\\
{\cal B}_K&=\left\{\omega : \SS_K\in\left(\Phi _*\left(\alpha \right),\Phi
_*\left(\beta   \right)\right)\right\}.
\end{align*} 
As before, we have the consistency
\begin{align*}
\bar\vartheta_{K }\longrightarrow \vartheta _0\qquad {\rm as} \qquad
K\rightarrow \infty 
\end{align*}
and 
\begin{align*}
\Ex_{\vartheta _0}\left|\bar\vartheta_{K }-\vartheta _0 \right|^2\leq \frac{C}{K}.
\end{align*}
We need the equation for $\dot m\left(\vartheta ,t\right)$ and expression for
Fisher information $${\rm I}\left(\vartheta_0 \right)=\sigma
^{-2}a^2\Ex_{\vartheta _0}\dot m\left(\vartheta_0 ,t\right)^2 $$ in this
case. The filtering equations in the stationary regime are
\begin{align*}
{\rm d}m\left(\vartheta  ,t\right)
&=-\left[ f +\frac{ \gamma_* \left( \vartheta \right)a^2
  }{\sigma^2}\right] m\left(\vartheta ,t\right){\rm d}t + \frac{\gamma_*
  \left(\vartheta \right)a}{\sigma^2} {\rm d}X_t,\qquad m\left(\vartheta
,0\right)= \xi ,\\ 
\gamma _*\left(\vartheta \right)&=\frac{f
  \sigma ^2}{a^2}\left(\sqrt{1+\frac{\vartheta^2a^2}{f ^2\sigma
    ^2}}-1\right),\qquad\qquad  \xi  \sim {\cal N}\left(0,\gamma _*\left(\vartheta
\right)\right) .
\end{align*}
Therefore
\begin{align*}
{\rm d}\dot m\left(\vartheta  ,t\right)
&=-\left[ f +\frac{ \gamma_* \left( \vartheta \right)a^2
  }{\sigma^2}\right]\dot m\left(\vartheta ,t\right){\rm d}t  + \frac{\dot\gamma_*
  \left(\vartheta \right)a}{\sigma^2} \left[{\rm d}X_t-a  m\left(\vartheta
  ,t\right){\rm d}t\right]. 
\end{align*}
For $\vartheta =\vartheta _0$ 
\begin{align*}
{\rm d}m\left(\vartheta_0 ,t\right) &=- f m\left(\vartheta_0 ,t\right){\rm d}t
+ \frac{\gamma_* \left(\vartheta_0 \right)a}{\sigma} {\rm d}\bar W_t,\qquad
m\left(\vartheta ,0\right) \sim {\cal N}\left(0,\gamma _*\left(\vartheta_0
\right)\right),\\
 {\rm d}\dot m\left(\vartheta_0 ,t\right) &=-A \left(
\vartheta_0 \right)\dot m\left(\vartheta_0 ,t\right){\rm d}t +
\frac{\vartheta_0a}{\sigma A \left( \vartheta_0 \right)} {\rm d}\bar W_t,
\; \dot m\left(\vartheta_0 ,0\right) \sim{\cal N}\left(0,q\left(\vartheta_0
\right)\right),
\end{align*}
where
\begin{align*}
 A\left(\vartheta _0\right)=f +\frac{ \gamma_* \left( \vartheta_0 \right)a^2
 }{\sigma^2}=\sqrt{f^2+\frac{\vartheta _0^2a^2}{\sigma ^2}} , \qquad
 q\left(\vartheta_0 \right)=\frac{\vartheta _0^2a^2}{\sigma ^2A\left(\vartheta
   _0\right)^3}
\end{align*}
Since 
\begin{align*}
\dot m\left(\vartheta_0  ,t\right)=\dot m\left(\vartheta_0
,0\right)e^{-At}+\int_{0}^{t} e^{-A\left(t-s\right)} \frac{\dot\gamma_*
  \left(\vartheta_0\right)a}{\sigma} {\rm d}\bar W_s
\end{align*}
we obtain 
\begin{align*}
\Ex_{\vartheta _0}\dot m\left(\vartheta_0  ,t\right)^2=
\frac{\vartheta_0^2a^2}{2\sigma^2 A\left(\vartheta _0\right)^3}
\end{align*}
Therefore the Fisher information is
\begin{align*}
{\rm I}\left(\vartheta \right)= \frac{\vartheta_0^2a^4}{2\sigma^4
  A\left(\vartheta _0\right)^3}  . 
\end{align*}

Now we can write the One-step MLE-process  $\vartheta _t^\star, T^\delta
<t\leq T$ as follows 
\begin{align}
\label{21}
\vartheta _t^\star=\bar\vartheta _{T^\delta }+\frac{a}{ \sigma ^2
  \left(t-T^\delta \right) {\rm I}\left(\bar\vartheta _{T^\delta } \right)}
\int_{T^\delta }^{t}\dot m(\bar\vartheta _{T^\delta },s) \left[{\rm
    d}X_s-am(\bar\vartheta _{T^\delta },s) {\rm d}s\right].
\end{align}
If we  change the variables $t=\tau T$ and denote $\vartheta _{\tau
  T}^\star=\vartheta_T ^\star\left(\tau \right), T^{\delta -1}<\tau \leq 1 $,
then we obtain the same assertions  as in the Theorem \ref{T1}:
\begin{proposition}
\label{P2}
  One-step MLE-process $\vartheta_T ^\star=\left(\vartheta_T ^\star\left(\tau
  \right), T^{\delta -1}<\tau \leq 1 \right) $ with $\delta \in \left(1/2,
  1\right)$ is consistent: for any $\nu >0$ and any $\tau \in (0,1]$
\begin{align*}
\lim_{T\rightarrow \infty }-\Pb_{\vartheta _0}\left(\left|\vartheta_T
^\star\left(\tau \right)-\vartheta _0 
\right|>\nu \right)=  0, 
\end{align*}
and asymptotically normal
\begin{align*}
\sqrt{\tau T}\left(\vartheta_T ^\star\left(\tau \right)-\vartheta _0
\right)\Longrightarrow {\cal N}\left(0, {\rm I}\left(\vartheta_0 \right)^{-
  1}\right).
\end{align*}
\end{proposition}
{\bf Proof.} Similarly to \eqref{16}, we have exactly the same representation for the
estimator $\vartheta _t^\star $ as in \eqref{14}, with the only difference in the forms
of $\dot m\left(\vartheta ,t\right)$ and ${\rm I}\left(\vartheta \right)$. Thus the previous proof 
works in this case as well. 

It is possible to write the system of recurrent equations as in
\eqref{16}-\eqref{18}.

\section{One-Step MLE-process. Case $\vartheta =a$.}

It is clear that the suggested estimation approach also works for the partially observed system 
\begin{align}
\label{22}
{\rm d}X_t&=\vartheta  Y_t\,{\rm d}t+
\sigma\, {\rm d}W_t,\qquad\qquad X_0=0, \\ {\rm
  d}Y_t&=-f Y_t\,{\rm d}t+b \,
           {\rm d}V_t,\qquad\qquad Y_0=\xi
\label{23}
\end{align}
where the unknown parameter is the drift $\vartheta=a$.

 The function
\begin{align*}
\hat \Phi\left(\vartheta \right)=\frac{b^2\vartheta
  ^2}{f^3}\left[e^{-f}-1+f\right] +\sigma ^2,\qquad \alpha <\vartheta <\beta 
\end{align*}
is  strictly increasing and the corresponding preliminary estimator
$\bar\vartheta _K$ admits the same asymptotic properties as in the preceding section. 

The filtering  equations are
\begin{align*}
{\rm d}m\left(\vartheta  ,t\right)
&=-\left[ f +\frac{\hat  \gamma \left( \vartheta \right)\vartheta^2
  }{\sigma^2}\right] m\left(\vartheta ,t\right){\rm d}t + \frac{\hat \gamma
  \left(\vartheta \right)\vartheta}{\sigma^2} {\rm d}X_t,\qquad m\left(\vartheta
,0\right)= \xi ,\\ 
\hat \gamma \left(\vartheta \right)&=\frac{f
  \sigma ^2}{\vartheta ^2}\left(\sqrt{1+\frac{\vartheta^2b^2}{f ^2\sigma
    ^2}}-1\right),\qquad\qquad  \xi  \sim {\cal N}\left(0,\hat \gamma\left(\vartheta
\right)\right) .
\end{align*}
Therefore
\begin{align*}
{\rm d}\dot m\left(\vartheta  ,t\right)
&=-\left[ f +\frac{\hat  \gamma \left( \vartheta \right)\vartheta^2
  }{\sigma^2}\right]\dot m\left(\vartheta ,t\right){\rm d}t  + \frac{\hat{\dot\gamma}
  \left(\vartheta \right)\vartheta +\hat{\gamma}
  \left(\vartheta \right)}{\sigma^2} {\rm d}X_t\\
&\qquad - \frac{\left[\hat{\dot\gamma}
  \left(\vartheta \right)\vartheta^2 +2\hat{\gamma}
  \left(\vartheta \right)\vartheta \right]}{\sigma^2} m\left(\vartheta ,t\right){\rm d}t
\end{align*}
To calculate Fisher information ${\rm I}\left(\vartheta
_0\right)=\sigma ^{-2}\Ex_{\vartheta _0}\left[ m\left(\vartheta_0
  ,t\right)+\vartheta _0\dot m\left(\vartheta_0  ,t\right)\right]^2$  we write
the representations 
\begin{align*}
m\left(\vartheta _0,t\right)&=\frac{\hat \gamma \left(\vartheta _0
  \right)\vartheta _0}{\sigma}\int_{0}^{t}e^{-f\left(t-s\right)}{\rm d}\bar
W_s +o\left(1\right), \qquad A=f +\frac{\hat \gamma \left( \vartheta_0
  \right)\vartheta_0^2 }{\sigma^2},\\ \dot m\left(\vartheta_0
,t\right)&=\frac{\hat{\dot\gamma} \left(\vartheta_0 \right)\vartheta_0
  +\hat{\gamma} \left(\vartheta_0
  \right)}{\sigma}\int_{0}^{t}e^{-A\left(t-s\right)}{\rm d}\bar W_s\\ &\qquad
\qquad -\frac{\hat{\gamma} \left(\vartheta_0 \right)\vartheta
  _0}{\sigma^2}\int_{0}^{t}e^{-A\left(t-s\right)}m\left(\vartheta
_0,s\right){\rm d}s+o\left(1\right) .
\end{align*}
In the  last integral we change the order of integration
\begin{align*}
&\int_{0}^{t}e^{-A\left(t-s\right)}m\left(\vartheta _0,s\right){\rm
    d}s=\frac{\hat \gamma \left(\vartheta _0 \right)\vartheta _0}{\sigma}
  \int_{0}^{t}e^{-A\left(t-s\right)}\int_{0}^{s}e^{-f\left(s-r\right)}{\rm
    d}\bar W_r\; {\rm d}s\\ &\qquad \qquad =\frac{\hat \gamma \left(\vartheta
    _0 \right)\vartheta _0}{\sigma}
  e^{-At}\int_{0}^{t}\left(\int_{r}^{t}e^{\left(A-f\right)s}{\rm d}s \right)
  e^{fr}{\rm d}\bar W_r \\ &\qquad \qquad =\frac{\hat \gamma \left(\vartheta
    _0 \right)\vartheta _0}{\sigma\left(A-f\right)}
  e^{-At}\int_{0}^{t}\left(e^{\left(A-f\right)t}-e^{\left(A-f\right)r}\right)
  e^{fr}{\rm d}\bar W_r \\ &\qquad \qquad =-\frac{\hat \gamma \left(\vartheta
    _0 \right)\vartheta _0}{\sigma\left(A-f\right)} \int_{0}^{t}
  e^{-A\left(t-r\right)}{\rm d}\bar W_r +o\left(1\right).
\end{align*}
Introduce notations
\begin{align*}
M\left(\vartheta _0\right)=\frac{\hat{\dot\gamma} \left(\vartheta_0
  \right)\vartheta_0 +\hat{\gamma} \left(\vartheta_0 \right)}{\sigma},\quad
N\left(\vartheta _0\right)= \frac{\hat{\gamma} \left(\vartheta_0
  \right)^2\vartheta _0^2}{\sigma^3\left(A-f\right)},\quad Q\left(\vartheta
_0\right)=\frac{\hat \gamma \left(\vartheta _0 \right)\vartheta _0}{\sigma} .
\end{align*}
Then $\dot m\left(\vartheta _0,t\right)$ we can write as follows
\begin{align*}
\dot m\left(\vartheta
_0,t\right)=\int_{0}^{t}e^{-A\left(t-s\right)}\left[M\left(\vartheta _0\right)
 +N\left(\vartheta _0\right)\right]{\rm d}\bar W_s+o\left(1\right).
\end{align*}
Hence
\begin{align*}
m\left(\vartheta
_0,t\right)+\vartheta _0\dot m\left(\vartheta
_0,t\right)&=Q\left(\vartheta _0\right)\int_{0}^{t}e^{-f\left(t-s\right)}
 {\rm d}\bar W_s\\
&\qquad +\left[M\left(\vartheta _0\right)
 +N\left(\vartheta _0\right)\right]\int_{0}^{t}e^{-A\left(t-s\right)}{\rm
   d}\bar W_s+o\left(1\right) 
\end{align*}
Therefore  the Fisher information in this problem is the function
\begin{align*}
{\rm I}\left(\vartheta _0\right)= \frac{Q\left(\vartheta _0\right)^2 }{2f} +
\frac{\left[M\left(\vartheta _0\right) +N\left(\vartheta
    _0\right)\right]^2}{2A}+\frac{2\left[M\left(\vartheta _0\right)
    +N\left(\vartheta _0\right)\right]Q\left(\vartheta
  _0\right)}{A+f}.
\end{align*}
Having the preliminary estimator $\bar\vartheta _{T^\delta }$, expression for
Fisher information ${\rm I}\left(\vartheta _0\right) $ and the equation for
$\dot m\left(\vartheta ,t\right)$ we can construct the One-step MLE-process
$\vartheta _t^\star, T^\delta <t\leq T $ of the same form 
as in \eqref{14}, with $a$ replaced by
$\bar\vartheta _{T^\delta }$.  

This estimator has the same asymptotic properties: it is consistent and
asymptotically normal
\begin{align*}
\sqrt{\tau T}\left(\vartheta _T^\star\left(\tau \right)-\vartheta
_0\right)\Longrightarrow {\cal N}\left(0,{\rm I}\left(\vartheta _0\right)^{-1}
\right).
\end{align*} 
The proof follows the same pattern as in the previous cases.  

\section{Discussion}

The results, presented above, can be developed in several directions by means of
already known approaches. 
\begin{enumerate}
\item It is interesting to find preliminary estimator in the cases of unknown
parameters $\vartheta =\left(f,b,a\right)$.  Of course, with one statistic
$\SS_{T^\delta }$ it is impossible and we need at least three different
statistics. 

Consider the case of two-dimensional parameter $\vartheta =\left(f,b\right)$
or $\vartheta =\left(f,a\right)$
and two statistics
\begin{align*}
\SS_K=\frac{1}{K}\sum_{k=1}^{K}\left[X_k-X_{k-1}\right]^2,\qquad
\RR_K=\frac{1}{K}\sum_{k=1}^{K}\left[X_k-X_{k-1}\right] \left[X_{k-1}-X_{k-2}\right] .
\end{align*}
The limits are 
\begin{align*}
\SS_K\longrightarrow \Phi \left(\vartheta
\right)=\frac{a^2b^2}{f^3}\left[e^{-f}-1+f\right],\qquad \RR _K\longrightarrow
\Xi\left(\vartheta \right)=\frac{a^2b^2}{2f^3}\left[e^{-f}-1\right]^2 .
\end{align*}
Therefore
\begin{align*}
\QQ_K=\frac{\SS_K}{ \RR _K}\longrightarrow \frac{2\left[e^{-f}-1+f\right]
}{\left[e^{-f}-1\right]^2}. 
\end{align*}
The function
\begin{align*}
\phi \left(x\right)=\frac{2\left[e^{-x}-1+x\right]
}{\left[e^{-x}-1\right]^2},\quad x> 0 
\end{align*}
is strictly increasing and $\lim_{x\rightarrow 0}\phi \left(x\right)=1 $,
$\lim_{x\rightarrow \infty }\phi \left(x\right)=\infty  $. Therefore, the
parameter $f$ can be estimated with the help of the statistic $\QQ_K$:
\begin{align*}
\QQ_K=\phi\left(f_K^*\right) .
\end{align*}
Having this estimator the second parameter, say, $a$ or $b$ can be
obtained as solution of one of these equations
\begin{align*}
\SS_K=\Phi \left(f_K^*,a_K^*\right),\qquad {\rm or}\qquad \SS_K=\Phi
\left(f_K^*,b_K^*\right), 
\end{align*}
with obvious notation. 
As soon as we have a consistent preliminary estimator, say, $\bar \vartheta
_{T^\delta}=\left(f_{T^\delta}^*,b_{T^\delta}^*\right)$ and explicit
expression for the information matrix ${\rm I}\left(\vartheta \right)$, then 
\begin{align*}
\vartheta _t^\star=\bar \vartheta _{T^\delta}+\left(t-T^\delta \right)^{-1}
          {\rm I}\left(\bar \vartheta _{T^\delta} \right)^{-1}\int_{T^\delta
          }^{t}\frac{a\,\dot m\left(\bar\vartheta _{T^\delta},s\right)}{\sigma
            ^2}\left[{\rm d}X_s- 
            m\left(\bar\vartheta _{T^\delta},s\right){\rm d}s\right] .
\end{align*}
Recall that such processes were studied in \cite{Kut17}.

\item The One-step MLE-process has {\it learning interval} $\left[0,T^\delta
  \right] $ with $\delta \in(\frac{1}{2},1]$. It can be interesting to have such
  process with {\it shorter} learning. This can be done with the help of
  another construction called Two-step MLE-process introduced in
  \cite{Kut17}. Let us recall this construction using the model of observation
  \eqref{7}-\eqref{8}. The first preliminary estimator $\bar \vartheta
  _{T^\delta }$ is constructed using the observations $X^{T^\delta
  }=\left(X_t,0\leq t\leq T^\delta \right)$ with $\delta \in (1/3,1/2]$
    (shorter learning interval). The second preliminary {\it
      estimator-process} $\vartheta _{t}^*, T^\delta <t\leq T $ is
\begin{align*}
\vartheta _{t}^*=\bar \vartheta _{T^\delta }+\frac{a}{\sigma ^2t{\rm I}\left(\bar
  \vartheta _{T^\delta} \right) }\int_{T^\delta }^{t}\dot
  m(\bar\vartheta _{T^\delta},s)\left[{\rm d}X_s-
  m(\bar\vartheta _{T^\delta},s){\rm d}s\right].
\end{align*}
The Two-step MLE-process is 
\begin{align*}
\vartheta _{t}^{\star\star}=\vartheta _{t}^*+\frac{a}{\sigma ^2t{\rm
    I}\left(\vartheta _{t}^* \right) }\int_{T^\delta }^{t}\dot m(\bar\vartheta
_{T^\delta},s)\left[{\rm d}X_s- m(\vartheta _{t}^*,s){\rm d}s\right].
\end{align*}
Following the same arguments as in the proof of Theorem 2 in \cite{Kut17} it
can be shown that
\begin{align*}
\sqrt{\tau T}\left(\vartheta _{T}^{\star\star}\left(\tau \right)-\vartheta _0
\right)\Longrightarrow {\cal N}\left(0,{\rm
    I}\left(\vartheta _0 \right)^{-1} \right),
\end{align*}
where $\vartheta _{T}^{\star\star}\left(\tau \right)=\vartheta _{\tau
  T}^{\star\star} $.

The learning interval $\left[0,T^{\delta }\right]$  can be made even shorter
if  $\delta \in (1/4,1/3]$. In this case we use Three-step MLE-process (see
  details in \cite{Kut17}). 

\item Consider the model \eqref{7}-\eqref{8} and the estimator-process
  $\vartheta _T^\star\left(\tau \right), \left[\kappa ,1\right]$, where
  $\kappa >0$. Let us denote by $ {\cal P}_T$ the measure induced by the process
  $$
\zeta _T\left(\tau \right)=\sqrt{T{\rm I}\left(\vartheta _0\right) }\left(
\vartheta _T^\star\left(\tau   \right)-\vartheta _0\right), \kappa \leq \tau \leq 1
$$ in the measurable space $\left({\scr C}\left[\kappa ,1\right], {\scr
  B}\right)$ of continuous on $\left[\kappa ,1\right]$ functions.  It is
possible to verify the weak convergence
\begin{align*}
{\cal P}_T\Longrightarrow {\cal P}
\end{align*}
where ${\cal P} $ corresponds to the Gaussian process $\zeta \left(\tau
\right),\left[\kappa ,1\right]$ with 
\begin{align*}
\Ex_{\vartheta _0}\zeta \left(\tau
\right)=0,\qquad \Ex_{\vartheta _0}\zeta \left(\tau_1\right)\zeta
\left(\tau_2\right)=\tau _1 \wedge \tau _2,
\end{align*}
i.e. $\zeta \left(\cdot  \right)$ is a Wiener process on the interval
$\left[\kappa ,1\right]$. 

The proof in similar situation can be found in \cite{Kut17}, Theorem 1. It
consists of proving convergence of the finite-dimensional
distributions 
\begin{align*}
\left(\zeta _T\left(\tau_1 \right),\ldots,\zeta_T \left(\tau_k
\right)\right)\Longrightarrow \left(\zeta \left(\tau_1 \right),\ldots,\zeta
\left(\tau_k \right)\right)
\end{align*}
and the estimate
\begin{align*}
\Ex_{\vartheta _0}\left|\zeta _T\left(\tau_1 \right)-\zeta_T \left(\tau_2
\right) \right|^4\leq C\left|\tau _2-\tau _1\right|^2,
\end{align*}
where the constant $C>0$ does not depend on $T$. The approach applied in the
present work allows us the direct verification these two conditions.

\end{enumerate}

{\bf Acknowledgment.} I would like to thank P. Chigansky for useful comments.

\end{document}